\documentclass[11pt,a4paper]{amsart}
\usepackage{amsmath,amsthm,amssymb}
\begin{document}

\newtheorem{theorem}{Theorem}[section]
\newtheorem{lemma}[theorem]{Lemma}
\newtheorem{proposition}[theorem]{Proposition}
\newtheorem{corollary}[theorem]{Corollary}

\theoremstyle{definition}
\newtheorem{definition}[theorem]{Definition}
\newtheorem{example}[theorem]{Example}
\newtheorem{formula}[theorem]{Formula}

\theoremstyle{remark}
\newtheorem{remark}[theorem]{Remark}

\renewcommand{\arraystretch}{1.2}

% definitons

\newcommand{\contr}{{\mspace{1mu}\lrcorner\mspace{1.5mu}}}
\newcommand{\per}{{\mspace{-1mu}\cdot\mspace{-1mu}}}
\newcommand{\de}{\partial}
\newcommand{\debar}{{\overline{\partial}}}
\newcommand{\desude}[2]{{\dfrac{\de #1}{\de #2}}}
\newcommand{\mapor}[1]{{\stackrel{#1}{\longrightarrow}}}
\newcommand{\ormap}[1]{{\stackrel{#1}{\longleftarrow}}}
\newcommand{\mapver}[1]{\Big\downarrow\vcenter{\rlap{$\scriptstyle#1$}}}
\newcommand{\bi}{\boldsymbol{i}}
\newcommand{\ide}[1]{{\mathfrak #1}}
\newcommand{\sA}{\mathcal{A}}
\newcommand{\Oh}{\mathcal{O}}
\newcommand{\sI}{\mathcal{I}}
\newcommand{\sF}{\mathcal{F}}
\newcommand{\sH}{\mathcal{H}}
\newcommand{\sL}{\mathcal{L}}
\newcommand{\DER}{{\mathcal{D}er}}

\newcommand{\C}{\mathbb{C}}
\newcommand{\Z}{\mathbb{Z}}
\newcommand{\N}{\mathbb{N}}

\newcommand{\MC}{\operatorname{MC}}
\newcommand{\Def}{\operatorname{Def}}
\newcommand{\hDef}{\widetilde{\operatorname{Def}}}
\newcommand{\Hom}{\operatorname{Hom}}
\newcommand{\Htp}{\operatorname{Htp}}
\newcommand{\End}{\operatorname{End}}
\newcommand{\Mor}{\operatorname{Mor}}
\newcommand{\Hilb}{\operatorname{Hilb}}
\newcommand{\Pic}{\operatorname{Pic}}
\newcommand{\Aut}{\operatorname{Aut}}
\newcommand{\image}{\operatorname{Im}}
\newcommand{\coker}{\operatorname{Coker}}

% end of definitions

\title{Lie description of higher
obstructions to deforming submanifolds}
\author{Marco Manetti}
\address{\newline Dipartimento di Matematica ``Guido Castelnuovo'',\hfill\newline
Universit\`a di Roma ``La Sapienza'',\hfill\newline
P.le Aldo Moro 5,
I-00185 Roma\\ Italy.}
\date{v2, \today}

\email{manetti@mat.uniroma1.it}
\urladdr{www.mat.uniroma1.it/people/manetti/}

\begin{abstract}
To every morphism $\chi\colon L\to M$ of differential graded Lie
algebras we associate  a functors of artin rings $\Def_\chi$ whose
tangent and obstruction spaces are respectively the first and
second cohomology group of the suspension of the mapping cone of $\chi$.\\
Such construction applies to Hilbert and Brill-Noether functors
and allow to prove with ease that every higher obstruction to
deforming a smooth submanifold of a K\"{a}hler
manifold is annihilated by the semiregularity map.
\end{abstract}

\maketitle

\section*{Introduction}

This paper is devoted, as its ancestors \cite{EDF} and \cite{CCK},
to develop algebraic tools that are useful to handle deformation
problems over a field of characteristic 0. The philosophy
underlying this project, due to Deligne,  Drinfeld and Kontsevich,
is that every "reasonable" deformation problem is the truncation
of an extended deformation problem which is governed by a
differential graded Lie algebra (DGLA) and every morphism of
deformation
theories is induced by an $L_{\infty}$-morphism.\\
Usually, the formal deformations of an algebro-geometric structure
are described by a functor of Artin rings \cite{Artinbook},
\cite{Sch}. The link between differential graded Lie algebras and
functors of Artin rings is given by the Maurer-Cartan equation.
More precisely, to every differential graded Lie algebra $L$ it is
associated the functor \cite{GoldMil1}, \cite{GoldMil2}, \cite{K},
\cite{ManettiDGLA}
\[\Def_L\colon \{\text{Local artinian rings}\}\to
\{\text{Sets}\}\]
\[ \Def_L(A)=\frac{\{x\in L^1\otimes \ide{m}_A\mid
dx+[x,x]/2=0\}}{\text{gauge  equivalence}}\]%
where $\ide{m}_A$ is the maximal ideal of $A$. The functor
$\Def_L$ depends only by the quasiisomorphism class of $L$.\\
In practice, the advantage of this approach concerns especially
the study of higher obstructions. Classically these obstructions
were studied using Massey products (see e.g. \cite{Douady},
\cite{palamodov}) but, as well explained in \cite{SchSta}, "Massey
product structures can be very helpful, though they are in general
described in a form that is unsatisfactory."\\
In most cases, the replacing of the Massey product structure with
the appropriate DGLA structure allow to
prove new results and easier proofs of old results.
On the other side it is not always easy to determine the right
DGLA governing a deformation problem.\\
The first goal of this paper is to develop some algebraic tools
that are useful to relate differential graded Lie algebras and
\emph{semitrivial deformations};  the paper \cite{ran4}
works on the same goal with a different approach.\\
The embedded deformations of a submanifold is the most classical
example of semitrivial deformation problem. Roughly speaking, we
have a manifold $X$ and a submanifold $Z\subset X$. The idea is to
consider the embedded deformation of $Z$ in $X$ as the
deformations of the inclusion map $Z\hookrightarrow X$ inducing a
trivial deformation of $X$. One other classical example is the
Brill-Noether functor, i.e. deformations of a bundle inducing
trivial deformation of its cohomology.\\
The key point is the definition of a functor $\Def_\chi$
associated to a morphism $\chi\colon L\to M$ of differential
graded Lie algebras and the verification that this functor has the
same formal properties of $\Def_L$; in particular the \emph{implicit
function theorem} (Theorem \ref{thm.basic}) still holds. Then,
using the trick of path-objects we see that there exists a
new DGLA $H$ such that $\Def_\chi=\Def_H$. This construction
allows for example to determine a differential graded Lie algebra
that governs embedded deformations and avoid the use, as in \cite{BC}, of higher
order differential operators.\\
The second goal of this paper is to check the utility of this
approach is a series of concrete examples; in particular we are
able to prove the following result (Theorem \ref{thm.main2}).

\begin{theorem}\label{thm.main}
Let $Z$ be a smooth closed submanifold of a compact K\"{a}hler
manifold $X$ and let $\omega$ be a closed differential
$(p,q)$-form on $X$ such that
$\omega_{|Z}=0$. Then the obstructions of embedded deformations of $Z$ inside $X$ are
contained in the kernel of the contraction map
\[ \contr\omega\colon H^1(Z,N_{Z|X})\to H^{q+1}(Z,\Omega_Z^{p-1}).\]
\end{theorem}

As applications of Theorem \ref{thm.main} we get
the vanishing of all higher obstructions to embedded
deformations under Bloch's semiregularity map  (Corollary
\ref{cor.semireg}) and the unobstructness of Lagrangian
submanifolds of a holomorphic symplectic
variety (Corollary \ref{cor.lagrangian}).\\

{\em Acknowledgements.}
I'm indebted with Domenico Fiorenza for several,
precious and useful discussions on these topics.
I  also  thank  Donatella Iacono for
having carefully read the  preliminary versions of this
paper.

\subsection*{General notation}
We always work over the filed $\C$, although most of the results
of algebraic nature are valid over an arbitrary field of
characteristic 0. Unless otherwise specified the symbol $\otimes$
denotes the tensor
product over $\C$.\\
The term  DGLA  means differential $\Z$-graded Lie algebras, while
dg-algebra means a differential
$\Z$-graded, graded-commutative and associative algebra.\\
Unless otherwise specified, every complex manifold is assumed compact and connected.
For every  complex manifold $X$ we denote by:
\begin{itemize}

\item $T_X$ the holomorphic tangent bundle on $X$.

\item $\sA_X^{p,q}$ the sheaf of differentiable $(p,q)$-forms of $X$. More generally
if $E$ is a holomorphic vector bundle on $X$ we denote by
$\sA_X^{p,q}(E)$  the sheaf of differentiable $(p,q)$-forms of $X$
with values in $E$ and by $A_X^{p,q}(E)=\Gamma(X,\sA_X^{p,q}(E))$ the space
of its global sections.

\item For every submanifold $Z\subset X$, we denote by $N_{Z|X}$ the normal bundle of $Z$ in $X$.

\end{itemize}

\section{Background}

Let $\mathbf{Set}$ be the category of sets and  $\mathbf{Art}$ the
category of local Artinian $\C$-algebras $(A, \ide{m})$ with
residue field
$A/\ide{m}=\C$.\\
Following \cite{Sch}, by a functor of Artin rings we intend a
covariant functor $\sF\colon \mathbf{Art}\to\mathbf{Set}$ such
that $\sF(\C)=\{\text{one element}\}$.\\
By the term Schlessinger's condition we mean one of the four
conditions $(H_{1}),\ldots,(H_{4})$ described in Theorem 2.1 of
\cite{Sch}.

The functors of Artin rings are used to describe infinitesimal
deformations of algebro-geometric structures.

\begin{example}\label{ex.hilbertfunctor}
Let $X$ be a complex manifold and let $Z\subset X$ be an analytic
subvariety defined by a sheaf of ideals $\sI\subset \Oh_X$.\\
The infinitesimal embedded deformations of $Z$ in $X$ are
described by the Hilbert functor $\Hilb^Z_X\colon\mathbf{Art}\to
\mathbf{Set}$,
\[ \Hilb_X^Z(A)=\{\text{ideal sheaves }\sI_A\subset
\Oh_X\otimes_{\C}A, \text{ flat over }A\text{ such that
}\sI_A\otimes_A\C=\sI\}.\]
\end{example}

Let $K=(\oplus K^i, d, [,])$ be a differential graded Lie algebra;
given $(A,\ide{m})\in\mathbf{Art}$, the set of Maurer-Cartan
elements with coefficients in $A$ is by definition
\[ \MC_{K}(A)=\left\{ x\in K^{1}\otimes \ide{m}\mid
dx+\frac{1}{2}[x,x]=0\right\},\]%
where the DGLA structure on $K\otimes \ide{m}$ is given by the
natural extension of the differential $d$ and the bracket $[~,~]$ defined as
\[ d(v\otimes a)=d(v)\otimes a,\qquad [v\otimes a,w\otimes b]=[v,w]\otimes ab.\]
Since $K^{0}\otimes \ide{m}$ is a nilpotent Lie algebra, its
exponential group $\exp(K^{0}\otimes \ide{m})$ can be defined as the set
$\{e^a\mid a\in K^{0}\otimes
\ide{m}\}$; we have $e^a e^b=e^{a\bullet b}$, where $\bullet$
is the Baker-Campbell-Hausdorff product (see e.g. \cite{Ja},
\cite{defomanifolds}).\\
The {\em gauge action} $\exp(K^{0}\otimes \ide{m})\times
\MC_{K}(A)\to \MC_{K}(A)$ is given explicitly by the formula
\[
e^a\ast w= w+\sum_{n\ge 0}\frac{[a,-]^{n}}{(n+1)!}([a,w]-da).
\]

\begin{remark}
The vector space $dK^{-1}$ is a Lie subalgebra of $K^0$ and
$\exp(dK^{-1}\otimes \ide{m})$ is contained in the stabilizer of
$0$. More generally, for every $x\in \MC_{K}(A)$ the subspace
$(d+[x,-])K^{-1}\otimes\ide{m}$ is a Lie subalgebra of
$K^0\otimes\ide{m}$ and the group
\[ S_x(A)=\{e^{[x,h]+dh}\mid h\in K^{-1}\otimes \ide{m}\}\]
is contained in the stabilizer of $x$. It is easy to verify that
for $a\in K^0\otimes \ide{m}$ we have $e^aS_x(A)e^{-a}=S_y(A)$,
where
$y=e^{a}\ast x$.\\
We shall call the group $S_x(A)$ the \emph{irrelevant stabilizer}
of $x$.
\end{remark}

The functor $\Def_{K}\colon \mathbf{Art}\to\mathbf{Set}$ is
defined as the quotient of Maurer-Cartan by  the gauge action,
\[\Def_{K}(A)=\frac{\MC_{K}(A)}{\exp(K^{0}\otimes \ide{m})}.\]%
The functor
$\Def_{K}$ satisfies the
Schlessinger's conditions $(H_{1})$, $(H_{2})$ (see
\cite{Sch},\cite{FM1}) and its tangent space
$\Def_{K}(\C[\epsilon]/(\epsilon^{2}))$
is naturally isomorphic to $H^{1}(K)$.\\
We point out for later use that, if $H^1(K)=0$, then
the functor $\Def_{K}$ is trivial and then for every
$x\in \MC_K(A)$ there exists $a\in K^0\otimes \ide{m}$ such that
$x=e^a\ast 0$.\\

\begin{definition}
We say that a functor of Artin rings $\sF\colon \mathbf{Art}\to
\mathbf{Set}$ is governed by a differential graded Lie algebra $K$
if $\sF$ is isomorphic to $\Def_K$.
\end{definition}

\begin{example}
Let $E$ be a holomorphic vector bundle on a complex manifold
$X$, then the functor of infinitesimal deformations of $E$ is governed by the
differential graded Lie algebra
\[K=\oplus_{i\ge 0}K^i,\qquad  K^i=A^{0,i}_X(\End(E)),\]
endowed with the Dolbeault differential and the natural bracket.
More precisely if $e,g$ are  local holomorphic sections of
$\End(E)$ and $\phi,\psi$  differential forms we define $d(\phi
e)=(\overline{\de}\phi)e$, $[\phi e,\psi g]=\phi\wedge\psi[e,g]$.
We refer to \cite{fuka}, \cite[Chap. VII]{Koba}, \cite[Sec.
9.4]{GoldMil1}, \cite[Pag. 238]{DK} for the proof that $\Def_K$ is
isomorphic to the  functor of infinitesimal deformations of $E$.
Here we only note that, for every $(A,\ide{m})\in\mathbf{Art}$ and
every $x\in \MC_K(A)$, the associated deformation of $E$ over
$Spec(A)$ is the bundle whose sheaf of holomorphic section is the
kernel of
\[ \debar+x\colon \sA^{0,0}_X(E)\otimes A\to \sA^{0,1}_X(E)\otimes A.\]
\end{example}

\bigskip

\section{Mapping cone of DGLA morphisms}

Let $\chi\colon L\to
M$ be a morphism of differential graded Lie algebras;
by definition the mapping cone of $\chi$ is the complex (see e.g. \cite{HartRes})
$\operatorname{Cone}(\chi)$, where
$\operatorname{Cone}(\chi)^i=L^{i+1}\oplus M^{i}$ and the differential is given by the
formula
\[ L^{i+1}\oplus M^i\ni (l,m)\mapsto (-dl,-\chi(l)+dm)\in L^{i+2}\oplus M^{i+1}.\]

We are interested to the \emph{suspension of the mapping cone} (SMC) of $\chi$;
it is  the differential graded vector space
$(C_\chi,\delta)$, where
$C_\chi^i=L^i\oplus M^{i-1}$ and the differential $\delta$ is defined as
\[ \delta(l,m)=(dl, \chi(l)-dm).\]
The projection $C_\chi\to L$  is a morphism of complexes and
there exist boundary operators $H^{i}(M)\to H^{i+1}(C_\chi)$ giving
a long exact sequence
\[\cdots\to H^i(C_\chi)\to H^i(L)\mapor{\chi}H^i(M)\to H^{i+1}(C_\chi)\to\cdots\]

In general does not exist any bracket on the suspended mapping cone
making $C_\chi$ a DGLA and  the projection $C_\chi\to L$ a Lie
morphism. Nevertheless, there exists a natural notion of
Maurer-Cartan equation and gauge action and we are able to define
the
associated deformation functor.\\
For every local Artinian $\C$-algebra $(A,\ide{m})$ we define
\[ \MC_{\chi}(A)=
\left\{(x,e^a)\in (L^1\otimes\ide{m})\times \exp(M^0\otimes
\ide{m})\mid dx+\frac{1}{2}[x,x]=0,\;e^a\ast\chi(x)=0\right\},\]
\[ \Def_\chi(A)=\frac{\MC_{\chi}(A)}{\exp(L^0\otimes\ide{m})
\times \exp(dM^{-1}\otimes\ide{m})},\]%
where the gauge action is given by the formula
\[ (e^l, e^{dm})\ast (x,e^a)=(e^l\ast x, e^{dm}e^ae^{-\chi(l)})=(e^l\ast x,
e^{dm\bullet a\bullet (-\chi(l))}).\]%
We now analyze the main properties of the functors $\Def_\chi$.

\medskip
{\sc Functoriality.} Every commutative diagram of differential
graded Lie algebras
\begin{equation}\label{eqn.square}
\begin{array}{ccc}
L&\mapor{f}&H\\
\mapver{\chi}&&\mapver{\eta}\\
M&\mapor{f'}&I\end{array}
\end{equation}%
induces a natural transformation of functors
$\Def_{\chi}\to \Def_{\eta}$
and a morphism of complexes of vector spaces
$C_\chi\to C_\eta$.

\begin{theorem}[Inverse function theorem]\label{thm.basic}
If the diagram {\rm (\ref{eqn.square})}
induces a quasi-isomorphism $C_\chi\to C_\eta$,
then the natural transformation $\Def_{\chi}\to
\Def_{\eta}$ is an isomorphism.
\end{theorem}

The proof is quite long and, for the clarity of exposition, it  is postponed  at the end of
Section \ref{sec.ESD}.\\

If $K$ is the kernel of a morphism $\chi\colon L\to
M$ of DGLA, then there exist natural transformations
\[ \Def_K\to \Def_{\chi}\mapor{\pi} \Def_{L}.\]
If $\chi$ is surjective, then $\Def_K\simeq\Def_{\chi}$; if $M$ is
acyclic, then $\Def_{\chi}\simeq \Def_{L}$.\\

\medskip
{\sc Tangent space.} The tangent space of a functor $\sF$ is by
definition the space $\sF(\C[\epsilon])$, where $\epsilon^2=0$.\\
Therefore
\[ \MC_\chi(\C[\epsilon])=
\left\{(x,e^a)\in (L^1\otimes\C\epsilon)\times
\exp(M^0\otimes\C\epsilon)\mid
dx=0,\;e^a\ast\chi(x)=\chi(x)-da=0\right\}\]%
\[ \simeq \left\{(x,a)\in
L^1\times M^0\mid dx=0,\;\chi(x)-da=0\right\}=
\ker(\delta\colon C_\chi^1\to C_\chi^2).\]%
Two elements $(x,a),(y,b)\in \ker\delta$ are gauge equivalent if
and only if there exists $(c,z)\in L^0\times M^{-1}$ such that
\[ y=x-dc,\qquad b=dz+a-\chi(c),\quad
 \text{or equivalently}\quad (x,a)-(y,b)=\delta(c,z).\]%
In conclusion the tangent space of $\Def_\chi$ is isomorphic to
$H^1(C_\chi)$.

\medskip
{\sc Obstructions.} The obstruction space of $\Def_\chi$ is
naturally contained in $H^2(C_\chi)$. Consider in fact a small
extension in $\mathbf{Art}$
\[0\mapor{} \C\mapor{}
A\mapor{\alpha}B\mapor{}0\]%
and let  $(x,e^q)\in \MC_\chi(B)$.\\
Since $\alpha$ is surjective there exists a pair  $(y,e^p)\in
L^1\otimes{\ide{m}_A}\times \exp(M^0\otimes \ide{m}_A)$ such that
$\alpha(y)=x$ and $\alpha(p)=q$.\\
Setting
\[ h=dy+\frac{1}{2}[y,y]\in L^2\otimes\C,\qquad
r=e^p\ast \chi(y)\in M^1\otimes\C\]%
we have $\delta(h,r)=0$. In fact,
\[dh=\frac{1}{2}d[y,y]=[dy,y]=[h,y]-\frac{1}{2}[[y,y],y].\]
By Jacobi identity $[[y,y],y]=0$, while $[h,y]=0$ because the
maximal ideal of $A$ annihilates the kernel of $\alpha$; therefore
$dh=0$.\\
Since $\chi(y)=e^{-p}\ast r= r+e^{-p}\ast 0$, we have
\[\chi(h)=d(r+e^{-p}\ast 0)+\frac{[r+e^{-p}\ast 0,r+e^{-p}\ast
0]}{2}=\]%
\[=dr+d(e^{-p}\ast 0)+\frac{[e^{-p}\ast 0,e^{-p}\ast
0]}{2}=dr,\]%
where the last equality follows from the fact that
$e^{-p}\ast 0$
satisfies the Maurer-Cartan equation in $M\otimes\ide{m}_A$.\\
The cohomology class of $(h,r)$ in $H^2(C_\chi)$ is well defined
and is exactly the obstruction of lifting $(x,e^q)$ to a
$\MC_\chi(A)$.\\
A straightforward computation shows that the primary obstruction
is the quadratic map
\[ H^1(C_\chi)\to H^2(C_\chi),\qquad (x,a)\mapsto
\frac{1}{2}([x,x], [a,\chi(x)]),\]%

\begin{lemma}\label{lem.critsmooth}
Let $\chi\colon L\to M$ be a morphism of DGLA and assume that
either
\begin{enumerate}

\item $L^1=0$, or

\item $[L^1,L^1]=0$ and $[M_0,dM_0]=0$
\end{enumerate}
Then $\Def_\chi$ is smooth.
\end{lemma}

\begin{proof}
In the first case the Maurer-Cartan equation
$dl+[l,l]/2=\chi(l)-e^{-a}\ast 0=0$
reduces to $e^{-a}\ast 0=0$ which is equivalent to $da=0$.\\
In the second case the Maurer-Cartan equation
$dl+[l,l]/2=\chi(l)-e^{-a}\ast 0=0$ reduces to $dl=\chi(l)-da=0$
and then $\MC_\chi$ is smooth.
\end{proof}

\begin{proposition}\label{prop.vaniobs}
Let
\[\begin{array}{ccc}
L&\mapor{}&H\\
\mapver{\chi}&&\mapver{\eta}\\
M&\mapor{}&I\end{array}
\]
be a commutative diagram of morphisms of differential
graded Lie algebras.\\
If the functor $\Def_\eta$ is smooth, then the obstruction space of $\Def_{\chi}$ is
contained in the kernel of $H^2(C_\chi)\to H^2(C_\eta)$.
\end{proposition}

\begin{proof}
The horizontal arrows are morphisms of DGLA;  
the explicit description of obstructions given above implies that 
the two morphism $\Def_{\chi}\to\Def_{\eta}$ and $H^2(C_\chi)\to H^2(C_\eta)$
commute with obstruction maps.\\
The smoothness of $\Def_{\eta}$ means 
that every obstruction vanishes in $H^2(C_\eta)$.
\end{proof}

\bigskip
\section{Example: deformations of cohomology groups}

For every pair $V,W$ of graded vector spaces we denote by
\[ \Hom^{*}(V,W)=\bigoplus_{i\in\Z}\Hom^{i}(V,W),\]
where
\[\Hom^{i}(V,W)=\{ f\colon V\to W\mid f(V_j)\subset
W_{i+j}\}=\prod_{j}\Hom(V_{j},W_{i+j}).\]

If $(V,d_V)$ and $(W,d_W)$ are  differential graded vector spaces,
then $\Hom^{*}(V,W)$ has a natural differential
\[ \delta\colon \Hom^{i}(V,W)\to \Hom^{i+1}(V,W),\qquad
\delta(f)=d_Wf-(-1)^{\deg(f)} fd_V.\]

Every $f\in \Hom^{i}(V,W)$ such that $\delta(f)=0$ induces
naturally a morphism $H(f)\in \Hom^{i}(H^*(V),H^*(W))$ and  for
every $i$ the natural transformation
\[ H\colon H^i(\Hom^{*}(V,W))\mapor{}\Hom^{i}(H^*(V),H^*(W))\]
is an isomorphism.\\
Every pair of morphisms of complexes $p\colon V'\to V$, $q\colon
W\to W'$ induces  a morphism of differential graded vector spaces
\[ r\colon \Hom^{*}(V,W)\to \Hom^{*}(V',W').\]
Notice that if $p,q$ are quasiisomorphisms, then  $r$ is a
quasiisomorphism.

Given a differential graded vector space $(V,d)$, the spaces
\[ \Hom^{*}(V,V),
\qquad  \Hom^{+}(V)=\bigoplus_{i\ge 0}\Hom^{i}(V,V)\]%
are differential graded Lie algebras with bracket
 \[ [f,g]:=f\circ g-(-1)^{\deg(f)\deg(g)}g\circ f.\]
Note that $\delta(f)=[d,f]$ and $f\in \Hom^{1}(V,V)$ satisfies the
Maurer-Cartan equation if and only if $(d+f)^{2}=0$.

\begin{lemma}\label{lem.formality}
If $V,W$ are quasiisomorphic complexes of $\C$-vector spaces, then
$\Hom^{*}(V,V)$ and $\Hom^{*}(W,W)$ are quasiisomorphic as
differential graded Lie algebras.
\end{lemma}

\begin{proof} Since every complex of vector spaces contains its
cohomology as a subcomplex, it is not restrictive to assume
$V\subset W$ a subcomplex.\\
The subspace $K=\{f\in \Hom^{*}(W,W)\mid f(V)\subset V\}$ is a
differential graded Lie subalgebra and there exists a natural
morphism of DGLA $K\to \Hom^{*}(V,V)$. There exists a commutative
diagram of complexes with exact rows
\[ \begin{array}{ccccccccc}
0&\mapor{}&K&\mapor{\alpha}&\Hom^{*}(W,W)&\mapor{}&\Hom^{*}(V,W/V)&\mapor{}&0\\
&&\mapver{\beta}&&\mapver{\gamma}&&\mapver{Id}&&\\
0&\mapor{}&\Hom^{*}(V,V)&\mapor{}&\Hom^{*}(V,W)&\mapor{}&\Hom^{*}(V,W/V)&
\mapor{}&0\end{array}\]%
Since $\Hom^{*}(V,W/V)$ is acyclic and
$\gamma$ is a quasiisomorphism, it follows that also $\alpha$ and
$\beta$ are quasi-isomorphisms.\end{proof}

It is clear that the inclusion of DGLA's $\Hom^{+}(V,V)\subset
\Hom^{*}(V,V)$ induces an
isomorphism of deformation functors.\\
Denoting, for notational simplicity,  $L=\Hom^{+}(V,V)$, for every
local Artinian $\C$-algebra $(A,\ide{m})$, the  Maurer-Cartan
elements  $\MC_{L}(A)$ are exactly the deformations of  the
differential $d$ over $A$, while the group
$\exp(L^{0}\otimes\ide{m})$ is the group of automorphisms of the
graded $A$-module $V\otimes A$ lifting the identity on $V$.

\begin{definition} We shall say that
$x\in \MC_{L}(A)$ gives a deformation of $H^i(V,d)$ if the
cohomology group $H^{i}(V\otimes A,d+x)$ is  a flat $A$-module and the
projection onto the residue field induce an isomorphism
$H^{i}(V\otimes A,d+x)\otimes_{A}\C\simeq H^{i}(V,d)$.
\end{definition}

It is clear that the subset of $x\in \MC_{L}(A)$
giving a deformation of $H^i(V,d)$ is stable under the gauge
action. We recall that every flat module over an artinian ring is free.

\begin{lemma}\label{lem.flatness}
In the notation above, if $H^{i-1}(V\otimes A,d+x)$ and
$H^{i+1}(V\otimes A,d+x)$ are deformations of $H^{i-1}(V,d)$ and
$H^{i+1}(V,d)$ respectively, then also
$H^{i}(V\otimes A,d+x)$ is a deformation of $H^{i}(V,d)$.\end{lemma}

\begin{proof} This is standard (see e.g. \cite{Artinbook}).
\end{proof}

\begin{lemma}\label{lem.triviality}
In the notation above, assume that the complex $(V,d)$ is bounded, i.e.
\[ (V,d):\qquad 0\mapor{}V^{m}\mapor{d}V^{m+1}\mapor{d}
\cdots\mapor{d}V^{n}\mapor{}0\]
for some pair of integers $m\le n$. Then
$x\in MC_{L}(A)$  is gauge equivalent
to $0$ if and only if  the (entire) cohomology of $(V\otimes A, d+x)$
is a deformation of the cohomology of $(V, d)$.
\end{lemma}

\begin{proof}
First we note that every deformation of a vector space over $A$ is
trivial and then  the cohomology of $(V\otimes A, d+x)$
is a deformation of the cohomology of $(V, d)$
if and only if $H(V\otimes A,d+x)\simeq H(V\otimes A,d)=H(V,d)\otimes A$.\\
By definition, if $a\in L^{0}\otimes\ide{m}$, then
$e^{a}\ast x=e^{a}(d+x)e^{-a}-d$ and therefore
$x,y\in MC_{L}(A)$ are
gauge equivalent if and only if the  complex $(V\otimes A, d+x)$
is isomorphic to $(V\otimes A, d+y)$.\\
In particular if $x$ is gauge equivalent of $0$ then
the cohomology of $(V\otimes A, d+x)$ is isomorphic to the
cohomology of $(V\otimes A, d)$.\\
Conversely assume that $x\in MC_{L}(A)$ and
$H(V\otimes A,d+x)\simeq H(V\otimes A,d)$;
since the complex is bounded, the local flatness criterion
implies that  there exists an isomorphism of complexes
$(V\otimes A,d+x)\simeq (V\otimes A,d)$.\end{proof}

Choosing an index $i$, we consider the truncated complex
\[ (T_iV,d):\qquad V^{i-1}\mapor{d}V^{i}\mapor{d}
V^{i+1}.\]
Then there exists a natural morphism of graded vector spaces
\[ T_i\colon\Hom^{+}(V,V)\to \Hom^{*}(T_iV,T_iV),\]%
induced by the inclusion $T_iV\subset V$ and the projection $V\to
T_iV$. It is straightforward to check that $T_i$ is a morphism of
DGLA; moreover the commutative diagram
\[ \begin{array}{ccc}
\Hom^{+}(V,V)&\mapor{Id}&\Hom^{+}(V,V)\\
\mapver{T_i}&&\mapver{}\\
\Hom^{*}(T_iV,T_iV)&\mapor{}&0\end{array}\] induces a natural
transformation of functors $\Def_{T_i}\to \Def_L$.

\begin{proposition}
In the notation above, assume it is given a
deformation of the complex $x\in MC_L(A)$.
The class $[x]\in \Def_L(A)$ belongs
to the image of $\Def_{T_i}(A)\to \Def_L(A)$
if and only if the cohomology group
$H^i(V\otimes A, d+x)$ is a deformation of $H^i(V, d)$.
\end{proposition}

\begin{proof}
By definition, the image of $\Def_{T_i}(A)\to \Def_L(A)$ is
the kernel of
\[T_i\colon \Def_{L}(A)\to \Def_{\Hom^{*}(T_iV,T_iV)}(A).\]%
The proof is an immediate consequence of Lemma
\ref{lem.triviality} and Lemma \ref{lem.flatness}.
\end{proof}

\bigskip
\section{Example:  Brill-Noether functors of vector bundles}

Let $E,F$ be  holomorphic vector bundles on a compact complex
manifold $X$ of dimension $n$. We want to describe the
deformations $E_t$ of
$E$ such that $\dim H^i(F\otimes E_t)=\dim H^i(F\otimes E)$ for every index $i$.\\
We have seen that
the infinitesimal deformations of $E$ are governed by the
differential graded Lie algebra
\[L=\oplus_{i\ge 0}L^i,\qquad  L^i=A^{0,i}_X(\End(E)).\]
This means that for every $(B,\ide{n})\in \mathbf{Art}$
\[ \Def_L(B)\simeq\{ E_{B}\to X\times Spec(B)\mid
E_B\text{ is a deformation of }E\}/\sim\]
where $\sim$ denotes isomorphism of deformations.
Our attention is devoted to the subfunctor
$\sF\subset \Def_L$,
\[ \sF(B)=\{ E_{B}\to X\times Spec(B)\mid
H^i(F\otimes E_B)\text{ is a deformation of }H^i(F\otimes E)\}/\sim.\]

Every $x\in L^i$  induces naturally  morphisms of sheaves
\[ x\colon \sA^{0,j}_X(F\otimes E)\to \sA^{0,i+j}_X(F\otimes E).\]
Denoting by
\[ (A_X(F\otimes E),\debar):\quad 0\mapor{}A^{0,0}_X(F\otimes E)\mapor{\debar}
A^{0,1}_X(F\otimes E)\mapor{\debar}\cdots\mapor{\debar}
A^{0,n}_X(F\otimes E)\mapor{}0\]
the  Dolbeault complex of $F\otimes E$,
there exists a natural morphism of DGLA
\[ \chi\colon L\mapor{}\Hom^{*}(A(F\otimes E),A(F\otimes E)).\]

For every  $x\in \MC_L(B)$ we denote by
$E_{B,x}\to X\times Spec(B)$ the corresponding deformation of $E$.
We recall that
\[ \Oh(F\otimes E_{B,x})=\ker(\sA^{0,0}_X(F\otimes E)\otimes
B\xrightarrow{\debar+x}\sA^{0,1}_X(F\otimes E)\otimes B)\]
and then  the
complex of sheaves
\[  \sA^{0,0}_X(F\otimes E)\otimes B\xrightarrow{\debar+x}\sA^{0,1}_X(F\otimes E)\otimes B
\xrightarrow{\debar+x}\sA^{0,2}_X(F\otimes E)\otimes B\xrightarrow{\debar+x}\cdots
\]%
is a fine resolution of $F\otimes E_{B,x}$. As a consequence, the
cohomology of $F\otimes E_{B,x}$ is isomorphic to the cohomology
of the complex of free $B$-modules
\[
A^{0,0}_X(F\otimes E)\otimes B\xrightarrow{\debar+x}A^{0,1}_X(F\otimes E)\otimes B
\xrightarrow{\debar+x}A^{0,2}_X(F\otimes E)\otimes B\xrightarrow{\debar+x}\cdots
\]

\begin{lemma}
Consider the morphism of differential graded Lie algebras
\[ L:=A_X^{0,*}(\End(E)) \mapor{\chi}\Hom^{*}(A_X(F\otimes E),A_X(F\otimes E)).\]
Then the subfunctor $\sF$ is the image of
$\Def_\chi\to\Def_L$.\\
In particular, the tangent space of $\sF$
is the kernel of the natural map
\[ H^1(\End(E))\to \bigoplus_i\Hom(H^{i-1}(F\otimes E),H^{i}(F\otimes E)).\]
\end{lemma}

\begin{proof} The first part is an immediate consequence of  the results
of the previous sections.
The tangent space is the kernel of
\[ H^1(\End(E))\to H^1(\Hom^{*}(A_X(F\otimes E), A_X(F\otimes E))).\]
On the other hand, according to Lemma \ref{lem.formality} we have
\[H^1(\Hom^{*}(A_X(F\otimes E),A_X(F\otimes E)))=
\oplus_i\Hom(H^{i-1}(F\otimes E), H^{i}(F\otimes E)).\]%
\end{proof}

As an application we get a new proof of the following
smoothness theorem of Green and Lazarsfeld
\cite{GL2}.

\begin{theorem} In the notation above, if $X$ is
compact K\"{a}hler, $E\in \Pic^0(X)$ and $F$ is
a flat unitary bundle, then $\Def_\chi$ is smooth.
\end{theorem}

\begin{proof}
We first point out that, if $G$ is a flat unitary vector bundle on
$X$, then  it makes sense to consider the sheaves
$\overline{\Omega}^{i}_X(G)=\ker(\de)$ of $G$-valued
antiholomorphic differential forms and, by Hodge theory on $G$,
the inclusion of complexes
\[ (\Gamma(X,\overline{\Omega}^{*}_X(G)),0)\subset (A^{0,*}_X(G),\debar),\]
is an injective quasiisomorphism.\\
Since $E$ is flat unitary
it is not restrictive to assume $E=\Oh_X$;
the differential graded Lie algebra $L$ is then abelian
and isomorphic to the complex $(A^{0,*}_X,\debar)$.
Denoting by
\[M=\{a\in \Hom^*(A_X(F),A_X(F))\mid a(\Gamma(X,\overline{\Omega}^{*}_X(F)))\subset
\Gamma(X,\overline{\Omega}^{*}_X(F))\}\]%
we have a commutative diagram of morphisms of DGLA
\[
\begin{array}{ccccc}
L&\leftarrow&\Gamma(X,\overline{\Omega}^{*}_X)&
\xrightarrow{Id}&\Gamma(X,\overline{\Omega}^{*}_X)\\
\mapver{\chi}&&\mapver{}&&\mapver{\eta}\\
\Hom^*(A_X(F),A_X(F))&\leftarrow&M&\rightarrow&
\Hom^*(\Gamma(X,\overline{\Omega}^{*}_X(F)),\Gamma(X,\overline{\Omega}^{*}_X(F))).
\end{array}\]
The horizontal arrows are quasiisomorphisms and therefore
$\Def_\chi=\Def_\eta$ by Theorem \ref{thm.basic}.
Then, according to  Lemma \ref{lem.critsmooth}, the functor  $\Def_\eta$ is  smooth.
\end{proof}

\bigskip

\section{Example: the Hilbert functor of a smooth submanifold}
\label{sec.hilbert}

Let $X$ be a smooth complex manifold of dimension $n$ and denote
by $\DER^{p}(\sA^{0,*}_X,\sA^{0,*}_X)$  the  sheaf of
$\C$-derivations of degree $p$ of the sheaf of graded algebras
$(\sA^{0,*}_X,\wedge)$ (note that $\debar\in
\DER^{1}(\sA^{0,*}_X,\sA^{0,*}_X)$).\\
The DGLA structure on the sheaf $\DER^{*}(\sA^{0,*}_X,\sA^{0,*}_X)=\oplus_p
\DER^{p}(\sA^{0,*}_X,\sA^{0,*}_X)$ is induced by the standard bracket
\[ [f,g]:=f\circ g-(-1)^{\deg(f)\deg(g)}g\circ f,\]
and the differential
\[ df:=[\debar,f]=\debar\circ f-(-1)^{\deg(f)}f\circ\debar.\]

For every $p$ we interpret $\sA^{0,p}_X(T_X)$ as a subsheaf of
$\DER^{p}(\sA^{0,*}_X,\sA^{0,*}_X)$, where the inclusion is
described in local holomorphic coordinates $z_1,\ldots,z_n$ by
\[ (\phi\desude{~}{z_i})(fd\overline{z}_{j_1}\wedge\cdots\wedge
d\overline{z}_{j_k})=\desude{f}{z_i}\phi\wedge
d\overline{z}_{j_1}\wedge\cdots\wedge d\overline{z}_{j_k}.\]

We note that $\sA^{0,*}_X(T_X)=\oplus_p \sA^{0,p}_X(T_X)$ is a
sheaf of differential graded  Lie subalgebras of
$\DER^{*}(\sA^{0,*}_X,\sA^{0,*}_X)$. A straightforward computation
shows that, if $z_{1},\ldots,z_{n}$ are local holomorphic
coordinates, $I,J$ ordered subsets of $\{1,\ldots,n\}$,
$a=fd\overline{z}_{I}\desude{~}{z_{i}}$,
$b=gd\overline{z}_{J}\desude{~}{z_{j}}$, $f,g\in\sA^{0,0}_X$ then
\[ da=\debar f\wedge d\overline{z}_{I}\desude{~}{z_{i}},\qquad
\left[a, b\right]=d\overline{z}_{I}\wedge d\overline{z}_{J}
\left(f\desude{g}{z_{i}}\desude{~}{z_{j}}-
g\desude{f}{z_{j}}\desude{~}{z_{i}}\right).\]

Assume now that $i\colon Z\hookrightarrow X$ is the inclusion of
closed smooth complex
submanifold and denote by 
\[i^*\colon (\sA^{0,*}_X,\debar)\to
(\sA^{0,*}_Z,\debar)\] 
the morphism of sheaves of differential
graded algebras given by restriction of forms on $Z$.\\
We denote by 
\[\sA^{0,*}_X(T_X)(-\log Z)=\{ \eta\in \sA^{0,*}_X(T_X)\mid 
\eta(\ker(i^*))\subset \ker(i^*)\}.\]
We note that  $\sA^{0,*}_X(T_X)(-\log Z)$ is a sheaf of
differential graded Lie subalgebras of $\sA^{0,*}(T_X)$ and
there exists an exact sequence of fine sheaves 
\[ 0\to \sA^{0,*}_X(T_X)(-\log Z)\to \sA^{0,*}_X(T_X)\to \sA^{0,*}_Z(N_{Z|X})\to 0.\]

Let $(A,\ide{m})$ be a local Artinian $\C$-algebra, every $\eta\in
\sA^{0,0}_X(T_X)\otimes \ide{m}$ induces an automorphism
\[ e^\eta\colon \sA^{0,*}_X\otimes A\to \sA^{0,*}_X\otimes
A,\qquad e^\eta(h)=\sum_{n=0}^{\infty} \frac{\eta^n}{n!}(h).\]%
If $\eta\in \sA^{0,0}_X(T_X)(-\log Z)\otimes \ide{m}$, then
$e^{\eta}(\ker(i^*)\otimes A)=\ker(i^*)\otimes A$.

\begin{lemma}\label{lem.dgauge}
For every local Artinian $\C$-algebra $(A,\ide{m})$ and every
$\eta\in \sA^{0,0}(T_X)\otimes\ide{m}$ we have
\[ e^\eta\circ \debar\circ e^{-\eta}=\debar+e^\eta\ast 0\colon
\sA^{0,0}_X\otimes A\to \sA^{0,1}_X\otimes A,\]%
where $\ast$ is the gauge action on $\sA^{0,*}_X(T_X)\otimes\ide{m}$. In particular
\[ \ker(\debar+e^{\eta}\ast 0\colon
\sA^{0,0}_X\otimes A\to \sA^{0,1}_X\otimes A)=e^{\eta}(\Oh_X\otimes A).\]

\end{lemma}

\begin{proof} This follows from the definition of the gauge action
and the fact that $\sA^{0,*}_X(T_X)$  is a  subalgebra of
$\DER^{*}(\sA^{0,*}_X,\sA^{0,*}_X)$.
\end{proof}

We denote by $A^{0,*}_X(T_X)(-\log Z)$  the differential graded Lie algebras of
global sections of the sheaf $\sA^{0,*}_X(T_X)(-\log Z)$, while, 
according to our
general notation, we denote by  $A^{0,*}_X(T_X)$ the DGLA of
global sections of $\sA^{0,*}_X(T_X)$. The differential graded Lie
algebra $A^{0,*}_X(T_X)$ is called the \emph{Kodaira-Spencer}
algebra of $X$.\\
The natural inclusion $\chi\colon A^{0,*}_X(T_X)(-\log Z)\to A^{0,*}_X(T_X)$  is a
morphism of differential graded Lie algebras and its cokernel is
isomorphic to the Dolbeault complex of $N_{Z|X}$; in particular
for every $i\ge 0$
\[ H^i(Z,N_{Z|X})\simeq H^i(A^{0,*}_X(T_X)/A^{0,*}_X(T_X)(-\log Z))\simeq H^{i+1}(C_{\chi}).\]
In the sequel of this section, just to avoid heavy formulas, we denote by $L_{Z|X}$ the 
differential graded Lie algebra $A^{0,*}_X(T_X)(-\log Z)$.\\

Consider now  the associated
functor $\Def_{\chi}$. Since $\chi$ is injective we
have, for every local Artinian $\C$-algebra $(A,\ide{m})$,
\[\MC_{\chi}(A)=\{e^\eta\in \Aut_A(\sA^{0,0}_X\otimes A)\mid
\eta\in A^{0,0}_X(T_X)\otimes\ide{m},\; e^{-\eta}\ast 0\in
L^1_{Z|X}\otimes \ide{m}\}.\]%
Under this identification the gauge action becomes
\[ \exp(L^0_{Z|X}\otimes\ide{m})\times \MC_{\chi}(A)\to \MC_{\chi}(A),\qquad
(e^\mu,e^\eta)\mapsto  e^\eta \circ e^{-\mu},\]%
and then
\[ \Def_{\chi}(A)=\frac{\MC_{\chi}(A)}{\exp(L^0_{Z|X}\otimes\ide{m})}.\]

\begin{theorem} There exists an isomorphism of functors
$\theta \colon \Def_{\chi}\to \Hilb^Z_X$.
\end{theorem}

\begin{proof} This is implicitly proved in \cite[Sec. 2]{BC} and \cite{clemens} using
the theory of transversely holomorphic trivialization (in the
relative case). Here we sketch a different proof.\\
Denote by $\sI=\Oh_X\cap \ker(i^*)$ the holomorphic ideal sheaf of $Z$ and
define
\[ \theta\colon \Def_{\chi}(A)\to \{\text{ideal sheaves of }\Oh_X\otimes_{\C}A\},
\qquad\theta(e^{\eta})=(\Oh_X\otimes A)\cap e^{\eta}(\ker(i^*)\otimes A).\]%
The morphism $\theta$ is well defined because, if $\mu\in
L^0_{Z|X}\otimes \ide{m}$, then
\[\theta(e^\eta\circ e^{-\mu})=(\Oh_X\otimes A)\cap e^\eta(e^{-\mu}(\ker(i^*)\otimes
A))=(\Oh_X\otimes A)\cap e^{\eta}(\ker(i^*)\otimes A).\]%
Next, we need to prove that $\theta(e^{\eta})$ is flat over $A$
and $\theta(e^{\eta})\otimes_A\C=\sI$; clearly we can prove the
same
properties for the sheaf $e^{-\eta}(\theta(e^{\eta}))$.\\
According to Lemma \ref{lem.dgauge}
\[e^{-\eta}(\Oh_X\otimes A)=\ker(\debar+e^{-\eta}\ast 0\colon
\sA^{0,0}_X\otimes A\to \sA^{0,1}_X\otimes A),\]%
and then
\[e^{-\eta}(\theta(e^{\eta}))=(\ker(i^*)\otimes A)\cap
\ker(\debar+e^{-\eta}\ast 0).\]%
Since flatness is a local property it is not restrictive to assume
$X$ a Stein manifold,  $H^1(X,T_X)=0$ and $H^0(X,T_X)\to
H^0(Z,N_{Z|X})$ surjective. This implies that $H^1(L_{Z|X})=0$ and
then the functor $\Def_{L_{Z|X}}$ is trivial. In particular there
exists $\mu\in L^0_{Z|X}\otimes \ide{m}$ such that $e^{-\mu}\ast
0=e^{-\eta}\ast 0$ and therefore
\[e^{-\eta}(\theta(e^{\eta}))=(\ker(i^*)\otimes A)\cap
\ker(\debar+e^{-\eta}\ast 0)=e^{-\mu}(\theta(e^{\mu}))=e^{-\mu}(\sI\otimes A).\]%
This proves that $\theta\colon \Def_\chi\to \Hilb^Z_X$.
It is well known (see e.g. \cite{kollar}) that the functor $\Hilb^Z_X$ is prorepresentable,
its tangent space is $H^0(Z,N_{Z|X})$ and its obstructions are in $H^1(Z,N_{Z|X})$. Therefore,
in order to prove that $\theta$ is an isomorphism it is sufficient to prove
that it is bijective on tangent space and injective on obstruction
space. This is a straightforward computation and it is left to the reader.
\end{proof}

In analogy with rational homotopy theory (see also next Remark
\ref{rem.whit}), it is possible to define the Whitehead product
\[ [~,~]_{W}\colon H^i(Z,N_{Z|X})\times H^j(Z,N_{Z|X})
\to H^{i+j+1}(Z,N_{Z|X})\]%
in the following way. For every cohomology class $a\in
H^i(Z,N_{Z|X})$ we denote by $\tilde{a}\in A_X^{0,i}(T_X)$ a
differential form that lifts $a$. This means in particular that
$d\tilde{a}\in L^{i+1}_{Z|X}$ and $\tilde{a}$ is
defined up to elements of $dA_X^{0,i-1}(T_X)+L^i_{Z|X}$.\\
Given $a\in H^i(Z,N_{Z|X})$ and  $b\in H^j(Z,N_{Z|X})$, we define
$[a,b]_W\in H^{i+j+1}(Z,N_{Z|X})$ as the cohomology class of
$\dfrac{1}{2}([\tilde{a},d\tilde{b}]-(-1)^i[d\tilde{a},\tilde{b}])$.
It is easy to verify that $[~,~]_{W}$ is well defined and induces
a structure of
graded Lie algebra on the space $\oplus V^i$, $V^i=H^{i-1}(Z,N_{Z|X})$.\\

\begin{corollary}
The primary obstruction map of $\Hilb^Z_X$ is equal to
\[ H^0(Z,N_{Z|X})\to H^1(Z,N_{Z|X}),\qquad a\mapsto \frac{1}{2}[a,a]_{W}.\]
\end{corollary}

\begin{proof} This follows from the description of the primary
obstruction map of the functor $\Def_\chi$.\end{proof}

\bigskip
\section{Review of extended deformation functors}

Here we review the basic definitions and properties of extended deformation functors,
introduced in \cite{EDF}.
We denote by:\begin{itemize}

\item $\mathbf{C}$  the category of all nilpotent finite dimensional
dg-algebras over $\C$.

\item $\mathbf{C}_0$ the full subcategory of $\mathbf{C}$ whose
objects are the dg-algebras $A\in\mathbf{C}$ with trivial
multiplication, i.e. $A\per A=0$.
\end{itemize}

In other words an object in $\mathbf{C}$ is a finite dimensional
complex $A=\oplus A_{i}$ endowed with a structure of dg-algebra
such that $A\per A\cdots(n\text{ factors})\cdots A=0$ for $n>>0$.
Note that if $A=A_{0}$ is concentrated in degree 0, then $A\in
\mathbf{C}$ if and only if $A$ is the maximal ideal of a local
artinian $\C$-algebra
with residue field $\C$.\\

There is an obvious equivalence
between $\mathbf{C}_0$
and the category of finite dimensional complexes of $\C$-vector spaces.

If $A\in\mathbf{C}$ and $I\subset A$ is a differential ideal, then
also $I\in \mathbf{C}$ and the inclusion $I\to A$ is a morphism of
dg-algebras.

\begin{definition}\label{VIII.5.1}
A {\em small extension} in $\mathbf{C}$ is a short exact sequence of complexes
\begin{equation*}
\label{smex}
0\mapor{}I\mapor{}A\mapor{\alpha}B\mapor{}0
\end{equation*}
such that $\alpha$ is a morphism in $\mathbf{C}$ and $I$ is an ideal of
$A$ such that  $AI=0$; in addition it is called {\em acyclic}
if  $I$ is an acyclic complex, or equivalently if
$\alpha$ is a quasiisomorphism.
\end{definition}

\begin{definition}\label{def.predef}
A covariant functor $F\colon \mathbf{C} \to
\mathbf{Set}$ is called a {\em predeformation functor} if the following
conditions are satisfied:
\begin{itemize}
\item[\ref{def.predef}.1:] $F(0)=\{\text{one element}\}$.

\item[\ref{def.predef}.2:] For every $A,B\in\mathbf{C}$,
the natural map
\[ F(A\times B)\to F(A)\times F(B)\]
is bijective.

\item[\ref{def.predef}.3:]  For every surjective
morphism
$\alpha\colon A\to C$ in $\mathbf{C}$, with $C\in \mathbf{C}_0$
an acyclic complex,
the natural morphism
\[ F(\ker(\alpha))\to F(A)\]
is bijective.

\item[\ref{def.predef}.4:]  For every pair of
morphisms
$\alpha\colon A\to C$, $\beta\colon B\to C$ in
$\mathbf{C}$, with $\alpha$ surjective,  the natural  map
\[F(A\times_C B)\to F(A)\times_{F(C)}F(B)\]
is surjective

\item[\ref{def.predef}.5:]
For every acyclic small extension
\[ 0\mapor{} I\mapor{} A\mapor{}B\mapor{}0\]
the induced map $F(A)\to F(B)$ is surjective.
\end{itemize}
\end{definition}

\begin{definition}\label{VIII.5.4}
A covariant functor $F\colon \mathbf{C} \to
\mathbf{Set}$ is called a
\emph{deformation functor} if it is a predeformation functor and
$F(I)=0$ for every acyclic complex $I\in \mathbf{C}_0$.
\end{definition}

\begin{example}\label{ex.expo}
Let  $L$ be a differential graded Lie algebra
and $A\in \mathbf{C}$; then the tensor product $L\otimes{A}$ has
a natural structure of nilpotent DGLA with
\begin{align*}
(L\otimes{A})^i=&\oplus_{j\in\Z}L^j\otimes A_{i-j}\\
d(x\otimes a)=&dx\otimes a+(-1)^{\deg(x)}x\otimes da\\
[x\otimes a, y\otimes b]=&(-1)^{\deg(a)\deg(y)}[x,y]\otimes ab\\
\end{align*}
Every morphism of DGLA, $L\to N$ and every morphism $A\to B$ in
$\mathbf{C}$ give a natural commutative diagram of morphisms of
differential graded Lie algebras
\[\begin{array}{ccc}
L\otimes{A}&\mapor{}&N\otimes{A}\\
\mapver{}&&\mapver{}\\
L\otimes{B}&\mapor{}&N\otimes{B}\end{array}\] The (extended)
exponential functor $\widetilde{\exp}_{L}\colon \mathbf{C}\to
\mathbf{Set}$ is defined as
\[\widetilde{\exp}_{L}(A)=\exp(H^{0}(L\otimes A)).\]%
It is a deformation functor in the sense of \ref{VIII.5.4}.
\end{example}

\begin{example}\label{ex.defo}
The (extended) deformation governed by $L$ are by definition given by the functor
\[ \hDef_L\colon \mathbf{C}\to \mathbf{Set},\qquad
\hDef_L(A)=\frac{\{x\in (L\otimes A)^1\mid dx+[x,x]/2=0\}}
{\text{Gauge action of }\exp((L\otimes A)^0)}.\]
It is proved in \cite{EDF} that $\hDef_L$ is a deformation functor
in the sense of Definition \ref{VIII.5.4}.
\end{example}

For every predeformation functor $F$ and every $A\in \mathbf{C}_0$ there
exists a natural structure of vector space on $F(A)$, where the sum and
the scalar multiplication are described by the maps
\[
A\times A\mapor{+}A\quad \Rightarrow\quad F(A\times A)=F(A)\times
F(A)\mapor{+}F(A)\]
\[ s\in\C,\quad A\mapor{\cdot s}A\qquad
\Rightarrow\qquad F(A)\mapor{\cdot s}F(A)\]
If $A\to B$ is a morphism in $\mathbf{C}_0$,
then $F(A)\to F(B)$ is $\C$-linear. Similarly if $F\to G$ is a
natural transformations of predeformation functors, the
map $F(A)\to G(A)$ is $\C$-linear for every
$A\in \mathbf{C}_0$.\\

\begin{definition}\label{VIII.5.6}
Let $F$ be a deformation functor and denote
$TF[1]^{n}=T^{n+1}F=F(\C\epsilon)$, where
$\epsilon$ is an indeterminate of degree $-n\in\Z$ such
that $\epsilon^{2}=0$.
The graded vector space
$TF[1]=\bigoplus_{n\in\Z}TF[1]^n$ is called  the
{\em tangent space} of $F$.\\
A natural transformation $F\to G$ of deformation functors is
called a {\em quasiisomorphism} if induces an isomorphism on tangent
spaces, i.e. if $T^{n}F\to T^{n}G$ is bijective for every $n$.
\end{definition}

For example, if $L$ is a differential graded Lie algebra, then
$T^{n}\widetilde{\exp}_{L}=H^{n-1}(L)$ and $T^{n}\hDef_{L}=H^{n}(L)$.

\begin{theorem}[Inverse function theorem]\label{VIII.6.4}
A natural transformation of deformation functors
is an isomorphism if and only if it is a quasiisomorphism.
\end{theorem}

\begin{proof} See \cite[Cor. 3.2]{EDF} or
\cite[Cor. 5.72]{defomanifolds}.\end{proof}

\begin{remark}\label{rem.whit}
It is proved in \cite{EDF} that for every deformation
functor $F$ it is defined the Whitehead product
on its tangent space
\[ [~,~]_{W}\colon TF[1]^i\times TF[1]^j\to
TF[1]^{i+j+1}\] inducing a graded Lie algebra structure on $\oplus
T^nF$. As in the topological case (see e.g. \cite[p. 111]{tanre}),
the product $[~,~]_{W}$ measure  the obstruction to lifting a map
from a wedge (of spheres in topology and of dg-fat points in
deformation theory) to  a product. In the topological analogy, the
space $\pi_i(X)\otimes \mathbb{Q}$ corresponds  to  $TF[1]^{-i}$,
where $F=\Mor(-,X)$.

\end{remark}

\bigskip
\section{Extended deformations from SMC}
\label{sec.ESD}

The aim of this section is to interpret the functors $\widetilde{\exp}_{L}$ and
$\hDef_{L}$ as special cases of a suitable functor $\hDef_{\chi}$,
where $\chi$ is a morphism of differential graded Lie algebras.

\begin{definition}
Let $\chi\colon L\to M$ be a morphism of differential graded Lie algebras.
The (extended) Maurer-Cartan functor $\MC_{\chi}\colon\mathbf{C}\to\mathbf{Set}$
is defined as
\[ \MC_{\chi}(A)=
\left\{(x,e^a)\in (L\otimes A)^{1}\times \exp((M\otimes A)^0)\mid
dx+\frac{1}{2}[x,x]=0,\;e^a\ast\chi(x)=0\right\},\]
where $*$ is the gauge action of $\exp((M\otimes A)^0)$ on $\MC_{\chi}(A)$
\[e^a\ast w= w+\sum_{n\ge 0}\frac{[a,-]^{n}}{(n+1)!}([a,w]-da).\]
\end{definition}

\begin{lemma}\label{VIII.5.7}
$\MC_{\chi}$ is a predeformation functor.
\end{lemma}

\begin{proof}
It is evident that $\MC_\chi(0)=0$ and for every pair of morphisms
$\alpha\colon A\to C$, $\beta\colon B\to C$ in $\mathbf{C}$ we
have
\[ \MC_\chi(A\times_C B)=\MC_\chi(A)\times_{\MC_\chi(C)}\MC_\chi(B)\]
and then  $\MC_{\chi}$ satisfies properties \ref{def.predef}.1,
\ref{def.predef}.2, \ref{def.predef}.3 and \ref{def.predef}.4.\\
Let $0\mapor{} I\mapor{} A\mapor{\alpha}B\mapor{}0$ be an acyclic small
extension and $(x,e^q)\in \MC_\chi(B)$.\\
Since $\alpha$ is surjective there exists a pair  $(y,e^p)\in
(L\otimes{A})^1\times \exp((M\otimes A)^0)$ such that
$\alpha(y)=x$ and $\alpha(p)=q$. Setting
\[ h=dy+\frac{1}{2}[y,y]\in (L\otimes{I})^2\]
we have
\[dh=\frac{1}{2}d[y,y]=[dy,y]=[h,y]-\frac{1}{2}[[y,y],y].\]
By Jacobi identity $[[y,y],y]=0$, while $[h,y]=0$ because $AI=0$;
therefore $dh=0$ and, since the complex  $L\otimes{I}$ is acyclic
by K\"{u}nneth formula, there
exists $s\in (L\otimes{I})^1$ such that $ds=h$.\\
The element $\hat{x}=y-s$ lifts $x$
and satisfies the Maurer-Cartan equation.\\
The element $z=e^p\ast \chi(\hat{x})$ belongs to $M\otimes I$ and
satisfies Maurer-Cartan equation. This means that $dz=0$ and then,
since $M\otimes I$ is acyclic, there exists $r\in (M\otimes I)^0$
such that
\[ e^re^p\ast\chi(\hat{x})=e^r\ast z=z-dr=0\]
and then $(\hat{x}, e^re^p)\in \MC_\chi(A)$ is a lifting of
$(x,e^q)$.
\end{proof}

\begin{definition}
The functor $\hDef_{\chi}\colon\mathbf{C}\to \mathbf{Set}$
is the quotient $\hDef_{\chi}(A)=\MC_{\chi}(A)/\sim$, where
\[ (x,e^q)\sim (e^a\ast x, e^{db}e^qe^{-\chi(a)}),\qquad a\in (L\otimes A)^0,\quad
b\in (M\otimes A)^{-1}.\]
\end{definition}

In other words, $\hDef_{\chi}(A)$ is the set of orbits of the action
\[ (\exp((L\otimes A)^0)\times \exp(d(M\otimes A)^{-1}))
\times MC_\chi(A)\to MC_\chi(A)\]
\[ ((e^a,e^{db}),(x,e^q))\mapsto (e^a\ast x, e^{db}e^qe^{-\chi(a)}).\]
Notice that $\exp(d(M\otimes A)^{-1}))$ is the irrelevant
stabilizer of $0$. We have:
\begin{itemize}
\item If $M=0$, then  $\hDef_{\chi}=\hDef_L$.
\item If $L=0$, then $\hDef_{\chi}=\widetilde{\exp}_M$.
\end{itemize}

\begin{theorem}\label{VIII.5.8}
$\hDef_{\chi}\colon \mathbf{C}\to \mathbf{Set}$
is a deformation functor with $T^i\hDef_{\chi}=H^i(C_\chi)$.
\end{theorem}

\begin{proof}
If $C\in \mathbf{C}_0$
then $L\otimes{C}$ and $M\otimes{C}$ are abelian differential graded Lie algebras and
\[ \MC_{\chi}(C)=\{(x,e^m)\mid dx=0,\; \chi(x)-dm=0\},\]
\[ \hDef_{\chi}(C)=\frac{\{(x,e^m)\mid dx=0,\; \chi(x)-dm=0\}}
{\{(-dy,e^{db-\chi(y)})\mid y\in (L\otimes C)^0,\; b\in (M\otimes C)^{-1}\}}.\]
Therefore $\hDef_{\chi}(C)$ is isomorphic  to the $H^1$ of the suspended mapping cone  of
$\chi\colon L\otimes C\to M\otimes C$. If $C$ is
acyclic then $\hDef_{\chi}(C)=0$; we have proved that $\hDef_{\chi}$ satisfies
the second condition of Definition \ref{VIII.5.4} and remains to prove that $\hDef_{\chi}$
is a predeformation functor.\\
Since $\hDef_{\chi}$ is the quotient of the predeformation functor $\MC_{\chi}$,
the conditions \ref{def.predef}.1  and \ref{def.predef}.5
are trivially verified.\\
It is clear from the definition that
$\hDef_{\chi}(A\times B)=\hDef_{\chi}(A)\times \hDef_{\chi}(B)$.\\
\emph{Verification of \ref{def.predef}.4}: Let $\alpha\colon A\to
C$, $\beta\colon B\to C$ morphism in $\mathbf{C}$ with $\alpha$ surjective.\\
Assume there are given $(a,e^p)\in \MC_{\chi}(A)$, $(b,e^q)\in
\MC_{\chi}(B)$ such that $\alpha(a,e^p)$ and $\beta(b,e^q)$ give
the same element in $\hDef_{\chi}(C)$; then there exist $u\in
(L\otimes{C})^0$ and $k\in (M\otimes C)^{-1}$ such that
\[\beta(b)=e^u\ast \alpha(a),\qquad
\beta(e^q)=e^{dk}\alpha(e^p)e^{-\chi(u)}.\]%
Let $v\in (L\otimes{A})^0$ be a lifting of $u$ and $h\in
(M\otimes{A})^{-1}$ be a lifting of $k$.\\
Replacing  $(a,e^p)$ with its gauge equivalent element $(e^v\ast
a, e^{dh}e^pe^{-\chi(v)})$, we may suppose
$\alpha(a,e^p)=\beta(b,e^q)$ and then the pair $((a,e^p),(b,e^q))$
lifts to $\MC_{\chi}(A\times_C B)$: this proves that the map
\[ \hDef_{\chi}(A\times_C B)\to \hDef_{\chi}(A)\times_{\hDef_{\chi}(C)}\hDef_{\chi}(B)\]
is surjective.\\
\emph{Verification of \ref{def.predef}.3} Assume
$C\in \mathbf{C}_0$ acyclic,
$\alpha\colon A\to C$ surjective and denote $D=\ker(\alpha)$.\\
Let $(a_1,e^{p_1}), (a_2,e^{p_2})\in MC_\chi(D)$, $u\in
(L\otimes{A})^0$ and $k\in (M\otimes A)^{-1}$ be
such that $a_2=e^u\ast a_1$ and $e^{p_2}=e^{dk}e^{p_1}e^{-\chi(u)}$.\\
We want to prove that there exist $v\in (L\otimes D)^0$
and $f\in (M\otimes D)^{-1}$ such that
$e^v\ast a_1=a_2$ and $e^{p_2}=e^{df}e^{p_1}e^{-\chi(v)}$.\\
Since $\alpha(a_1)=\alpha(a_2)=0$ and $L\otimes{C}$ is an abelian DGLA we
have $0=e^{\alpha(u)}\ast 0=0-d\alpha(u)$ and then $d\alpha(u)=0$.
Since $L\otimes C$ is acyclic
there exists $h\in
(L\otimes{A})^{-1}$ such that $d\alpha(h)=-\alpha(u)$ and $u+dh\in
(L\otimes{D})^0$.\\
Setting  $w=[a_1,h]+dh$, then $e^w$ belongs to the irrelevant stabilizer of $a_1$ and
therefore  $(e^u e^w)\ast a_1=e^u\ast a_1=a_2$. Writing $e^u e^w=e^v$,
we claim that $v\in L\otimes{D}$: in fact
\[v=u\bullet w\equiv u+w\equiv
u+dh\pmod{[L\otimes{A}, L\otimes{A}]}\] and since $A\per A\subset
D$ we have
$v=u\bullet w\equiv u+dh\equiv 0 \pmod{L\otimes{D}}$.\\
On the other hand, $e^{\chi(w)}$ belongs to the irrelevant
stabilizer of $\chi(a_1)$ and then there exists  $l\in (M\otimes
A)^{-1}$ such that $e^{dk}e^{p_1}e^{\chi(w)}=e^{dl}e^{p_1}$; we
can write
\[ e^{p_2}=e^{dk}e^{p_1}e^{-\chi(u)}=e^{dk}e^{p_1}e^{\chi(w)}e^{-\chi(w)}e^{-\chi(u)}=
e^{dl}e^{p_1}e^{-\chi(v)}.\]%
Since $\exp((M\otimes D)^0)$ is a subgroup of $\exp((M\otimes
A)^0)$
we have $e^{dl}\in \exp((M\otimes D)^0)$ and then $dl\in M\otimes D$.\\
The inclusion $M\otimes D\to M\otimes A$ is a quasiisomorphism and
then the cohomology class of $dl$ is trivial in $M\otimes D$.
There exists $f\in (M\otimes D)^{-1}$ such that $df=dl$ and then
$e^{p_2}=e^{df}e^{p_1}e^{-\chi(v)}$.
\end{proof}

It is clear that every commutative diagram of morphisms of
differential graded Lie algebras
\[ \begin{array}{ccc}
L&\mapor{f}&H\\
\mapver{\chi}&&\mapver{\eta}\\
M&\mapor{f'}&I\end{array}\]
induces a natural transformation of functors
\[ \hDef_{\chi}\to \hDef_{\eta}.\]
The inverse function theorem implies that such a natural
transformation is an isomorphism if and only if the pair $(f,f')$
induce a quasiisomorphism between the SMC of $\chi$
and $\eta$.

\begin{example}
Assume that $\chi\colon L\to M$ is a surjective morphism of
differential graded Lie algebras, then $N=\ker(\chi)$ is a
differential graded Lie algebra and the commutative diagram
\[ \begin{array}{ccc}
N&\hookrightarrow&L\\
\mapver{}&&\mapver{\chi}\\
0&\mapor{}&M\end{array}\]
induce an isomorphism $H^*(N)\mapor{\simeq} H^*(C_\chi)$.
Therefore the natural transformation $\hDef_N\to \hDef_{\chi}$ is an isomorphism.
\end{example}

The next theorem shows that, even if $\chi$ is not surjective, there exists a
differential graded Lie algebra $H$ such that $\hDef_H\simeq\hDef_{\chi}$.\\
Denote by $M[t,dt]=M\otimes \C[t,dt]$, where $\C[t,dt]$ is the
polynomial De Rham algebra of the affine line. More precisely
$\C[t,dt]=\C[t]\oplus \C[t]dt$,
$t$ has degree $0$, $dt$ has degree 1 and $d(p(t)+q(t)dt)=p'(t)dt$.\\
The inclusion $\C\to \C[t,dt]$ is a quasiisomorphism and then, by
K\"{u}nneth formula, also the inclusion $i\colon M\to M[t,dt]$ is
a
quasiisomorphism of differential graded Lie algebras.\\
Define, fore every $a\in \C$ the evaluation morphism
\[ e_a\colon M[t,dt]\to M,\qquad e_a(\sum m_it^i+n_it^idt)=\sum m_ia^i,\]
is a morphism of DGLA which is a left inverse of the inclusion $i$;
in particular $e_a$ is a surjective quasiisomorphism for every $a$.

\begin{theorem}\label{thm.atom}
For every morphism $\chi\colon L\to M$ of differential graded Lie algebras,
the fiber product
\[ H=\{(l,m)\in L\times M[t,dt]\mid e_0(m)=0,\, e_1(m)=\chi(l)\}\]
is a differential graded Lie algebra and $\hDef_H=\hDef_{\chi}$.\end{theorem}

\begin{proof} Setting
\[ K=\{(l,m)\in L\times M[t,dt]\mid  e_1(m)=\chi(l)\}\]
we have  a commutative diagram

\[ \begin{array}{ccc}
L&\mapor{f}&K\\
\mapver{\chi}&&\mapver{e_0}\\
M&\mapor{Id}&M\end{array}\qquad\qquad f(l)=(l,\chi(l)).\]
Passing to  SMC we get an isomorphism
$(f,Id)\colon H^*(C_\chi)\to H^*(C_{e_0})$ and $H$ is the kernel of the surjective
morphism $e_0$.
\end{proof}

\begin{proof}[Proof of Theorem \ref{thm.basic}]
By definition, for every $(A,\ide{m})\in\mathbf{Art}$ and every
$\chi\colon L\to M$ we have $\Def_\chi(A)=\hDef_\chi(\ide{m})$.
The proof follows from Theorem \ref{VIII.5.8} and inverse function
theorem \ref{VIII.6.4}.
\end{proof}

\bigskip
\section{A new look to Cartan formulas}

For every differential graded vector space $(V,d_V)$ and every
integer $i\in \Z$, we define the shifted differential graded
vector space $(V[i],d_{V[i]})$ by setting
\[ V[i]^j=V^{i+j},\qquad d_{V[i]}=(-1)^id_V.\]
Next, for every pair $V,W$ of differential graded vector spaces we
define
\[\Htp(V,W)=\Hom^*(V[1],W).\]%
In other terms, for every integer $i$ we have
\[ \Htp^i(V,W)=\Hom^i(V[1],W)=\Hom^{i-1}(V,W)\]
and the differential of $\Htp(V,W)$ is given by the formula
\[  \Htp^i(V,W)\ni f\mapsto\delta(f)=d_W f-
(-1)^i f d_{V[1]}=d_W f+ (-1)^i fd_V.\]%

Let $X$ be a complex manifold with Kodaira-Spencer algebra
$A^{0,*}_X(T_X)$
and denote by $(A_X,d)$ the De Rham complex of $X$, 
i.e. $A_X=\oplus_{p,q}A^{p,q}$ and $d=\de+\debar$.\\ 
The contraction map is the linear map
\[\bi\colon A^{0,*}_X(T_X)\to\Hom^{*}(A_X, A_X)\]
defined as
\[\bi_a(\omega)=a\contr \omega,\qquad a\in A^{0,*}_X(T_X), \quad \omega\in
A_X.\]%

 Note that $\bi\colon
A^{0,i}_X(T_X)\to\oplus_{h,l} \Hom_{\C}(A_X^{h,l},
A_X^{h-1,l+i})\subset \Hom^{i-1}(A_X,A_X)$ and then $\bi$ has
degree $-1$.

Denoting by $[\; ,\;]$ the standard bracket in the DGLA
$\Hom^{*}(A_X, A_X)$, we have the \emph{Cartan formulas} (for a
proof see e.g. \cite{CCK}, \cite{defomanifolds})
\[\bi_{da}=[\debar, \bi_a],\quad \bi_{[a,b]}=
[\bi_a,[\de,\bi_b]]=[[\bi_a,\de],\bi_b],\quad[\bi_a,\bi_b]=0.\]%
In order to interpret $\bi$ as a morphism of differential graded
Lie algebras we need to consider the DGLA given by the
differential graded vector space
$\Htp\left(\ker(\de),\dfrac{A_X}{\de A_X}\right)$ (whose
differential is $\delta(f)=\debar f+(-1)^{\deg(f)} f\debar$) with
the  bracket
\[  \{f,g\}=f\de
g-(-1)^{\deg(f)\deg(g)} g\de f.\]

\begin{proposition}
The linear map
\[ \bi\colon A_X^{0,*}(T_X)\to
\Htp\left(\ker(\de),\dfrac{A_X}{\de A_X}\right)\]%
is a morphism of differential graded Lie
algebras.\end{proposition}

\begin{proof} Straightforward consequence of Cartan
formulas.\end{proof}

Consider now a smooth closed submanifold $Z\subset X$ and denote
by $I_Z\subset A_X$ the graded subspace  of differential forms vanishing on $Z$.
Since
\[A^{0,*}_X(T_X)(-\log Z)\subset \{a\in
A^{0,*}_X(T_X)\mid \bi_a(I_Z)\subset I_Z\},\]%
we have a commutative diagram of morphisms of DGLA
\[ \begin{array}{ccl}
A^{0,*}_X(T_X)(-\log Z)&\mapor{\bi}&\left\{f\in
\Htp\left(\ker(\de),\dfrac{A_X}{\de A_X}\right)\mid
f(I_Z\cap\ker(\de))\subset \dfrac{I_Z}{I_Z\cap\de A_X}\right\}
\\
\mapver{\chi}&&\qquad\mapver{\eta}\\
A^{0,*}_X(T_X)&\mapor{\bi}& \Htp\left(\ker(\de),\dfrac{A_X}{\de
A_X}\right)\end{array}
\]

Notice that
\[ \coker(\eta)=
\Htp\left(I_Z\cap \ker(\de),\dfrac{A_Z}{\de A_Z}\right).\]%

\begin{lemma}\label{lem.ottodue} If the differential graded vector
spaces
$(\de A_X,\debar)$ and $(\de A_Z,\debar)$ are acyclic,
 then the functor $\Def_{\eta}$ is unobstructed.
 In particular the obstructions of $\Def_\chi=\Hilb^Z_X$ are
contained in the kernel of
\[
H^1(N_{Z|X})=H^2(C_{\chi})\mapor{\bi}H^2(C_{\eta})=\mathop{\oplus}_{i}
\Hom\left(H^{i}(I_Z\cap\ker(\de)),H^{i}\left(\dfrac{A_Z}{ \de
A_Z}\right)\right).\]

\end{lemma}

\begin{proof} We first note that the exact sequence
\[ 0\mapor{}I_Z\cap\de A_X\mapor{}\de A_X\mapor{}\de
A_Z\mapor{}0\]%
implies that also the complex $I_Z\cap\de A_X$ is acyclic.\\
For simplicity of notation denote
\[ K=\left\{f\in
\Htp\left(\ker(\de),\dfrac{A_X}{\de A_X}\right)\mid
f(I_Z\cap\ker(\de))\subset \dfrac{I_Z}{I_Z\cap\de A_X}\right\}.\]

The projection $\ker(\de)\to \ker(\de)/\de A_X$ induces a
commutative diagram
\[ \begin{array}{ccc}
\{f\in K\mid f(\de A_X)=0\}&\mapor{\alpha}&K\qquad
\\
\mapver{\mu}&&\mapver{\eta}\qquad\\
\Htp\left(\dfrac{\ker(\de)}{\de A_X},\dfrac{A_X}{\de
A_X}\right)&\mapor{\beta}& \Htp\left(\ker(\de),\dfrac{A_X}{\de
A_X}\right)\end{array}\]%

Since $\de A_X$ is acyclic, $\beta$ is a
quasiisomorphism of DGLA.\\
 Since
\[ \coker(\alpha)=\left\{f\in
\Htp\left(\de A_X,\dfrac{A_X}{\de A_X}\right)\mid f(I_Z\cap\de
A_X)\subset \dfrac{I_Z}{I_Z\cap\de A_X}\right\}.\]%
there exists an exact sequence
\[ 0\to \Htp\left(\dfrac{\de A_X}{I_Z\cap \de A_X},\dfrac{A_X}{\de A_X}\right)
\to\coker{\alpha}\to \Htp\left(I_Z\cap \de A_X,\dfrac{I_Z}{I_Z\cap
\de A_X}\right)\to 0.\]
Since the complexes $\dfrac{\de A_X}{I_Z\cap \de A_X}=\de A_Z$ and
$I_Z\cap \de A_X$ are both acyclic, also $\coker(\alpha)$ is acyclic
and then  $\alpha$ is a quasiisomorphism.
According to Theorem \ref{thm.basic} there exists an isomorphism of functors
$\Def_{\eta}=\Def_{\mu}$.\\
On the other side, both algebras on the first column are abelian
and then, by Lemma \ref{lem.critsmooth} the functor $\Def_{\mu}$
is smooth. The vanishing of obstructions follows from Proposition
\ref{prop.vaniobs}.
\end{proof}

\bigskip
\section{Examples and applications}
\label{sec.examples}

In the notation of previous section, the contraction operator gives a
morphism of complexes
\[ \bi\colon A_X^{0,*}(T_X)\to
\Htp\left(\ker(\de),A_X\right),\]%
where the differential on $\Htp\left(\ker(\de),A_X\right)$ is
$f\mapsto \debar f\pm f\debar$.
The inclusion $I_Z\cap \ker(\de)\subset \ker(\de)$ and the projection
$A_X\to A_Z$ give a morphism of complexes
\[ \pi_Z\colon \Htp\left(\ker(\de),A_X\right)\to
\Htp\left(I_Z\cap \ker(\de),A_Z\right).\]
The subalgebra $A^{0,*}_X(T_X)(-\log Z)$ is contained in the kernel of the composition
and therefore there exists a quotient map
\[ \pi_Z\bi\colon A_Z^{0,*}(N_{Z|X})\to \Htp\left(I_Z\cap \ker(\de),A_Z\right).\]
If $\omega\in I_Z$ is a closed $(p,q)$-form, then $\de \omega=\debar\omega=0$
and then gives a morphism of complexes
\[ \contr\omega\colon (A_Z^{0,*}(N_{Z|X}),\debar)\to (A_Z^{p-1,q+*},\debar),\qquad
\eta\contr\omega= \pi_Z\bi_{\eta}(\omega).\]

Therefore the Theorem \ref{thm.main} is completely equivalent to

\begin{theorem}\label{thm.main2}
In the above notation, if $X$ is compact K\"{a}hler, then the obstruction
to $\Hilb^Z_X$ are contained in the kernel of
\[ \pi_Z\bi\colon H^1(A_Z^{0,*}(N_{Z|X}))\to
H^1(\Htp\left(I_Z\cap \ker(\de),A_Z\right)).\]
\end{theorem}

\begin{proof}
Since $X$ is compact K\"{a}hler, the subcomplex
$\image(\de)=\de A_X\subset A_X$ is acyclic. In fact by $\de\debar$-lemma we have
\[\ker(\debar)\cap \image(\de)=\ker(\de)\cap \image(\debar)=\image(\de\debar)\]
and the equality
$\ker(\debar)\cap \image(\de)=\image(\de\debar)$ implies immediately that
$H^*_{\debar}(\image(\de))=0$.
Moreover, also  $Z$ is K\"{a}hler and then  the same conclusion
holds for the complex $\de A_Z$. In particular the projection
\[ \Htp\left(I_Z\cap \ker(\de),A_Z\right)\to
\Htp\left(I_Z\cap \ker(\de),\frac{A_Z}{\de A_Z}\right)\]
is a quasiisomorphism.
By Lemma \ref{lem.ottodue}, the obstruction space of $\Hilb_X^Z$
is  contained in the kernel of the linear map

\[\begin{array}{ccl}
H^1(N_{Z|X})&\mapor{\bi}&
H^1\left(\Htp\left(I_Z\cap \ker(\de),\dfrac{A_Z}{\de A_Z}\right)\right)\\
&&\qquad\qquad\Vert\\
&&\mathop{\oplus}_{i}
\Hom\left(H^{i}(I_Z\cap\ker(\de)),
H^{i}\left(\dfrac{A_Z}{\de A_Z}\right)\right)
\end{array}\]

\end{proof}

The interplay between  semiregularity and embedded deformations 
has been studied  by Severi \cite{severi} for curves on surfaces, 
by Kodaira and Spencer \cite{kodaspen} for submanifolds of codimension 1 and 
by S. Bloch \cite{bloch} for every  submanifolds of a projective variety.
They proved that 
if the semiregularity map is injective, then the corresponding 
embedded deformations are unobstructed.

\begin{corollary}\label{cor.semireg}
Let $Z$ be a smooth closed submanifold of codimension $p$ of a
compact K\"{a}hler manifold $X$. Then the obstruction space of
$\Hilb^Z_X$ is contained in the kernel of the semiregularity map
\[ H^1(Z,N_{Z|X})\to H^{p+1}(X,\Omega_X^{p-1}).\]
\end{corollary}

\begin{proof} Under our
assumption the semiregularity map 
can be defined in the following way: let $n$ be the
dimension of $X$ and denote by $\sH$ the space of harmonic forms
on $X$ of type $(n-p+1, n-p-1)$. By Dolbeault theorem and Serre
duality, the dual of $\sH$ is isomorphic to
$H^{p+1}(X,\Omega_X^{p-1})$.\\
The composition of the contraction map and integration on $Z$
gives a bilinear map
\[ H^1(Z,N_{Z|X})\times \sH\to H^{n-p}(Z,\Omega_Z^{n-p})\to\C,\qquad
(\eta,\omega)=\int_{Z}\eta\contr\omega\]%
which induces the semiregularity map
\[ H^1(Z,N_{Z|X})\to \sH^{\vee}=H^{p+1}(X,\Omega_X^{p-1}).\]
Since $\sH\subset I_Z\cap\ker(\de)\cap\ker(\debar)$, the proof
follows immediately from Theorem \ref{thm.main}.
\end{proof}

\begin{remark} It is not clear to me if Corollary \ref{cor.semireg}
is valid without the K\"{a}hler assumption; a closer look to the
proofs shows that the $\de\debar$-lemma on $Z$ is not required and
then  the K\"{a}hler assumption can be weakened to the validity of
the $\de\debar$-lemma on $X$. Therefore, according to \cite[Cor.
5.23]{DGMS}, the Corollary \ref{cor.semireg} holds for every
compact complex manifold $X$ which can be blown up to a K\"{a}hler
manifold; in particular it holds for every Moishezon
variety.\end{remark}

\begin{corollary}[Voisin \cite{voisin}]\label{cor.lagrangian}
Let $X$ be a K\"{a}hler and holomorphic symplectic variety and let
$Z\subset X$ be a holomorphic lagrangian submanifold. Then
$\Hilb_X^Z$ is smooth.
\end{corollary}

\begin{proof}
Let $\omega$ be the holomorphic symplectic $(2,0)$-form on $X$. If
$Z$ is a lagrangian submanifold then $\omega_{|Z}=0$ and the
contraction with $\omega$ gives an isomorphism of vector bundles
over $Z$
\[ \contr\omega\colon  N_{Z|X}\mapor{\simeq}\Omega^1_Z.\]
In particular the map
\[ \contr\omega\colon H^1(Z,N_{Z|X})\to H^{1}(Z,\Omega_Z^{1})\]
is injective and then, according to Theorem \ref{thm.main}, every
obstruction of $\Hilb_X^Z$ vanishes.\end{proof}

\end{document}